\documentstyle[11pt]{article}

\newtheorem{th}{Theorem}[section]
\newtheorem{lem}[th]{Lemma}
\newtheorem{pro}[th]{Proposition}
\newtheorem{rem}[th]{Remark}
\newtheorem{cor}[th]{Corollary}
\newtheorem{ex}[th]{Example}

 \newcommand{\ov}{\overline}

\newcommand{\uqp}{{\sf u}_q^+\!\left({\cal G}\right)}
\newcommand{\uqpsl}{{\sf u}_q^+\!\left({ sl_2}\right)}

\newcommand{\ot}{\otimes}
\newcommand{\oti}{{{\mbox{\tiny $\otimes$}}}}

\newcommand{\pr}{{\cal P}}
\newcommand{\si}{{\cal S}}

\newcommand{\zz}{{\bf Z}}
\newcommand{\cc}{{\bf C}}
\newcommand{\kpr}{K_0({\cal P}\Lambda)}
\newcommand{\kp}{K_0^p(\Lambda)}

\title{The projective class ring of half-quantum groups at roots of unity}

\author{Claude Cibils}

\begin{document} \date{} \maketitle

\begin{abstract}
We describe the structure of the Grothendieck ring of projective modules of
 basic Hopf algebras  using a positive integer determined by the composition
series of the principal indecomposable projective module.
\vskip2mm
\small
\noindent
{\bf 1991 Mathematics Subject Classification :} 16E20 and 16W30 \end{abstract}

\normalsize

\section{Introduction}

We consider finite dimensional Hopf  algebras which are basic and split,  which
 means that simple module are one-dimensional over the ground field. Some of
those algebras  have been called half-quantum groups and are useful in knot
theory and 3-manifold invariants, see \cite{crfr,ku}. Through the Drinfeld
quantum double (\cite{dr1}) they provide universal $R$-matrices and solutions to
the Yang-Baxter equation. Notice that the relation between the K-theory of a
Hopf algebra and the K-theory of its quantum double is not understood. Basic and
split Hopf algebras has been studied in \cite{hqg,acq,grso}.

The purpose of this  paper is to describe the  projective class ring of an
arbitrary basic and split finite  dimensional Hopf algebra, where the
multiplicative structure is given by the tensor product of modules obtained
through the comultiplication of the Hopf algebra.  We obtain for instance that
for $q$ a root of unity, the Grothendieck ring of $\uqpsl$ has nilpotent
elements. Notice that $\uqpsl$ has only a finite number of isomorphism classes
of indecomposable modules; by contrast Benson and Carlson proved in \cite{beca}
that a group algebra  verifying this property has no nilpotent elements in its
Grothendieck ring.

The tensor product of simple modules is often commutative, this happens to be
true for half-quantum groups at roots of unity.  With this assumption, the
structure of the complexification of the  projective class ring can be made more
precise than in the general case, see Theorem \ref{struc}.  The result relies on
a positive integer determined by the number of roots of unity where the
composition series of the projective cover of the trivial module do not vanish.
Its value for $sl_2$ is $1$, and the complexification of the projective class
ring of $sl_2$ is the product of two copies of the complex numbers with the
product of $n-1$ copies of the dual numbers.

The projective class ring of half-quantum groups at roots of unity are
commutative rings. Actually the entire Grothendieck ring  of $\uqpsl$ is
commutative, as a consequence of  direct computations made in \cite{aqqg}; there
is no theoretical result predicting this commutativity property since it is
known that half-quantum groups are not quasitriangular, and not even
quasi-cocommutative (see \cite{ra}).

\section{Basic and split Hopf algebras}

Let $\Lambda$ be a finite dimensional Hopf algebra over a field $k$, that is an
associative  algebra, provided with a comultiplication $\Delta:
\Lambda
\rightarrow \Lambda \otimes \Lambda$,
a counit $\varepsilon : \Lambda \rightarrow k$  and an antipode $S:\Lambda\rightarrow
\Lambda$  (see \cite{sw}).

The tensor product $M\otimes N$ over $k$ of two left $\Lambda$-modules $M$ and $N$ has  a
left  $\Lambda$-module structure  obtained through $\Delta$. The ground field
$k_\varepsilon$ equipped with the module structure given by $\varepsilon$ is the unit
element; since $M^*={\rm Hom}_k(M,k)$ is a right $\Lambda$-module, it turns into a left
module by composing the action with the antipode.

An associative algebra is split if the endomorphism rings of the simple modules are
 isomorphic to the ground field; this condition  is always satisfied if $k$ is
algebraically closed. The algebra is basic if each simple module has
multiplicity one in  its maximal semisimple quotient algebra. Each basic and
split Hopf algebra has a structure group attached, which corresponds to the
group of group-likes in the dual Hopf algebra.

\begin{lem}
Let ${\cal S}$ be the set of isomorphism classes of simple modules over a basic and split
finite dimensional Hopf algebra. Then $\cal S$ is a group.
\end{lem}
{\bf Proof:} Let $S$ and $T$ be simple modules. Since they are one-dimensional, their
tensor product is also one-dimensional, hence a simple module. The unit element in ${\cal
S}$ is the isomorphism  class of $k_{\varepsilon}$, while the inverse of $S$ is $S^*$.

\vskip3mm

\begin{ex} \rm
Let $G$ be a finite group and $k^G$ be the Hopf algebra of functions on it, that
is the dual of the group algebra $kG$. Its structure group  is $G$.
\end{ex}

\begin{ex}
\rm
Let $\uqp$ be the half-quantum group corresponding to a simple Lie algebra $\cal
G$ with Cartan matrix of size $t\times t$, and to a primitive root of unity of
order $n$. This basic and split algebra has an abelian structure group,
isomorphic to the product of $t$ copies of the cyclic group of order $n$.
\end{ex}

If the algebra is not basic, the tensor product of simple modules is not
semisimple in general. The structure group gives an easy proof of a particular
case of a result due to P. Schauenburg. Recall that associative algebras are
Morita equivalent if their module categories are equivalent. Hopf algebras are
Morita-monoidally equivalent if their module categories are equivalent as
monoidal categories, which means that there exist an equivalence of categories
compatible with tensor products, dual operations and the unit object. In
\cite{sc} Schauenburg shows that finite dimensional Hopf algebras
Morita-monoidally equivalent have isomorphic underlying associative algebras.

\begin{pro}
Let $\Lambda$ be a basic and split Hopf algebra, and let $\Lambda'$ be a Hopf
algebra Morita-monoidally equivalent to $\Lambda$. Then $\Lambda$ and $\Lambda'$
are isomorphic associative algebras.
\end{pro}
{\bf Proof:} Consider the structure group ${\cal S}$ of $\Lambda$ and the set
${\cal S'}$ of isomorphism classes of simple $\Lambda'$-modules. It follows from
the definition  that ${\cal S'}$ is also a group with respect to the tensor
product of modules, isomorphic to ${\cal S}$. Since the dimension over the
ground field is multiplicative on tensor products and since each simple module
has an inverse,  simple $\Lambda'$-modules are one-dimensional. Hence $\Lambda'$
is  basic and split algebra, and consequently isomorphic to $\Lambda$.
\vskip3mm

We recall now definitions which will be useful at the next section. The
Grothendieck group of an abelian category ${\cal C}$ is the free abelian group
with basis the isomorphism classes $[X]$ of objects $X$ of ${\cal C}$, divided
by the subgroup generated by elements $[X_2]-[X_1]-[X_3]$ provided by each split
exact sequence $0\rightarrow X_1 \rightarrow X_2 \rightarrow X_3
\rightarrow 0$ in ${\cal C}$. We denote the quotient group $K_0({\cal C})$, and
$K_0(\Lambda)$ denotes the Grothendieck group of the category of finitely
generated modules over a ring $\Lambda$. In case ${\cal C}$ is a monoidal
abelian category -- for instance the category of modules over  a Hopf algebra --
the free abelian group above is a ring. Tensoring a split exact sequence yields
to a split exact  sequence, the quoted subgroup is a two-sided ideal and the
Grothendieck group is a ring.

\begin{rem}
\rm
If $\Lambda$ is a finite dimensional algebra, $K_0(\Lambda)$ is a free abelian group with a basis
given by the isomorphism classes of indecomposable modules, as an immediate consequence of the
Krull-Schmidt Theorem (see for instance \cite{cure}).
\end{rem}

Consider the subgroup  generated by elements given by non necessarily split
exact sequences and denote the corresponding quotient group ${\overline
K}_0({\cal C})$. For modules over a finite dimensional algebra $\Lambda$, the
Jordan-Holder theorem insures that ${\overline K}_0({\Lambda})$ is free abelian
with a basis provided by the isomorphism classes of simple modules. Since
$\Lambda$  is a Hopf algebra over a field $k$, the functors $X\otimes -$ and
$-\otimes X$ are exact for any $\Lambda$ module $X$,  and ${\overline
K}_0(\Lambda)$ is a ring. If $\Lambda$ is a basic and split Hopf algebra, this
ring is clearly isomorphic to the group ring ${\bf Z}{\cal S}$, where ${\cal S}$
is the structure group of $\Lambda$.

There is a canonical ring surjection $K_0(\Lambda)\rightarrow {\overline
K}_0(\Lambda)$ called dimension vector which associates to a module $M$ the sum
${\overline M}$ of the quotients of one of its composition series. This map is
an isomorphism if and only if $\Lambda$ is a semisimple algebra.

\section{Projective class rings}

We  study $K_0^p(\Lambda)$, the subring of $K_0(\Lambda)$ generated by the
projective and semisimple modules. We recall the following standard results:

\begin{lem}
\label{projectives}
Let $\Lambda$ be a Hopf algebra, and let $P$ and $X$  be finitely generated
modules such that $P$ is projective. Then $P\otimes X$ is a projective module.
\end{lem}
{\bf Proof:} We  prove  that  $\Lambda \otimes X$ is free. The module
$\Lambda\otimes X$ is isomorphic to $\Lambda\otimes X_{\varepsilon}$, where
$X_\varepsilon$ is the vector space $X$ equipped with the trivial action given
by the counit $\varepsilon$. The result follows since $\Lambda\otimes
X_{\varepsilon}$ is free, of rank the dimension of $X$. An isomorphism of left
$\Lambda$-modules $\varphi:\Lambda\otimes X\rightarrow
\Lambda\otimes X_{\varepsilon}$ is given by $\varphi(a\oti x)=\sum a_{(1)}\oti
S(a_{(2)})x$, where $\Delta a = \sum a_{(1)}\oti a_{(2)}$.

\begin{pro}
Let $\Lambda$ be a finite dimensional basic and split Hopf algebra. The subring
$K_0^p(\Lambda)$ generated by semisimple and projective modules is a free
abelian group with basis the set of isomorphism classes of simple  and
indecomposable projective modules.
\end{pro}
{\bf Proof:} $K_0(\Lambda)$ is free abelian with basis the set of isomorphism
classes of indecomposable modules.  The described subgroup  is the smallest
containing semisimple and projective modules. Moreover this subgroup is an
ideal, due to Lemma \ref{projectives} and to the fact that isomorphism classes
of simple modules forms a group.
\vskip3mm
Let $\pr$ be  the set of isomorphism classes of indecomposable projective
modules. The ring ${\overline K}_0(\Lambda)$ is the integral group ring of the
structure group $\si$, and $\zz\pr$ is a $\zz\si$-bimodule.  It is well known
that $\si$ and $\pr$ are canonically bijective sets (see for instance
\cite{cure}) : each simple module $S$  has an unique projective cover $P_S$, a
projective module such that there exist an essential surjection
$\psi:P_S\rightarrow S$: for any strict submodule $X$ of $P_S$, the restriction
of $\psi$ to $X$ is not surjective.

\begin{lem}
\label{bimodule}
For a basic and split finite dimensional Hopf algebra, $\zz\pr$ is isomorphic to
$\zz\si$ as a $\zz\si$-bi\-mo\-du\-le.
\end{lem}
{\bf Proof:} Consider the above canonical bijection : we have to prove that
$T\ot P_S$ is isomorphic to $P_{T\ot S}$ for $T$ and $S$ simple modules. Since
$\psi:P_S\rightarrow S$ is an essential surjection, we obtain that $T\ot
P_S\rightarrow T\ot S$ is also essential since simple modules are
one-dimensional. Moreover,  $T\ot P_S$ is a projective module, and the unicity
of projective covers  insures that $T\ot P_S$ is the projective cover of $T\ot
S$.

\begin{lem}
Let $\Lambda$ be a finite dimensional basic and split Hopf  algebra which is not
semisimple.  Then $\kp$ is the direct sum of the subring $\zz\si$ and the
$\zz\si$-bimodule $\zz\pr$.

\end{lem}
{\bf Proof:} Since ${\overline K}_0(\Lambda)$ is the subring $\zz\si$,  we  need
to prove that $\si\cap\pr$ is empty; indeed  if there exist a simple and
projective module, we infer that any simple module is projective by Lemma
\ref{projectives} and the fact that $\si$ is a group. If all the simple modules
are projective modules, there are no non trivial extensions between them, and
the algebra is semisimple.

\vskip3mm

\begin{rem}
\rm
A semisimple basic and split Hopf algebras  is of the form $k^G$ for $G$ a
finite group. Such algebras are commutative and $\kp$ is the group algebra $\zz
G$ since the structure group of $k^G$ is $G$.

\end{rem}

The next purpose is to  describe the multiplicative structure of $\zz\pr$, in
other words to describe the $\zz\si$-bimodule map
 $$\zz\pr\ot_{\zz\si}\zz\pr\longrightarrow\zz\pr$$ induced by the tensor product
of modules over a basic and split Hopf algebra.
\begin{pro}
\label{devisse}
Let  $\Lambda$ be any finite dimensional Hopf algebra, and let $M$ and $P$ be
finitely generated left $\Lambda$-modules where $P$ is projective. Then $M\ot P$
is isomorphic to ${\ov M}\ot P$.
\end{pro}
{\bf Proof:} Recall that ${\ov M}$ is the semisimple module obtained by adding
the quotients of a composition series of $M$. It suffices to prove that for an
exact sequence $$0\rightarrow M_1\rightarrow M \rightarrow M_2\rightarrow 0$$ of
$\Lambda$-modules we have that $M\ot P$ is isomorphic to $M_1\ot P \ \oplus\
M_2\ot P$. By tensoring, we get an exact sequence $$0\rightarrow M_1\ot
P\rightarrow M\ot P
\rightarrow M_2\ot P\rightarrow 0$$ of $\Lambda$-modules. Since $M_2\ot P$ is
projective by Lemma \ref{projectives}, the sequence splits.
\vskip3mm
\begin{cor}
If the structure group of a basic and split finite dimensional Hopf algebra is
abelian, the  projective modules are central in $K_0(\Lambda)$.
\end{cor}
{\bf Proof:} We prove that any indecomposable projective module is in the center
of $K_0(\Lambda)$, more precisely we  have that $M\ot P_S$ and $P_S\ot M$ are
isomorphic modules for any finitely generated module $M$ and any simple module
$S$. From the proof of the preceding result we have that $P_S\ot M$ and $P_S\ot
{\ov M}$ are isomorphic modules. We need to prove that $P_S\ot T$ and $T\ot P_S$
are isomorphic for any simple module $T$. The module $P_S\ot T$ is the
projective cover of $S\ot T$ (see the proof of Lemma \ref{bimodule}) and $T\ot
P_S$ the projective cover of $T\ot S$. Since $\si$ is commutative, and
projective covers are unique up to isomorphism, we obtain the result.
\vskip3mm

Another way of formulating Proposition \ref{devisse} is to consider the
Gro\-then\-dieck group $\kpr$ of the category of projective finitely generated
$\Lambda$-modules. The natural two-sided action of $K_0(\Lambda)$ on $\kpr$
provides  an action of ${\ov K}_0(\Lambda)$ on $\kpr$.
\vskip3mm

\noindent
{\bf Definition.} Let $A$ be a ring; a multiplicative $A$-bimodule is an $A$-bimodule $M$ provided
with an associative $A$-bimodule map $\varphi :M\ot_A M\rightarrow M$.

\vskip3mm
If $M$ is a  multiplicative $A$-bimodule, then $A\oplus M$ has a natural ring
structure where $M$ is a two-sided ideal. Conversely, if $R$ is a ring whith a
subring $A$ and a decomposition $R=A\oplus I$ where $I$ is a two-sided ideal of
$R$, we obtain that $I$ is a multiplicative $A$-bimodule. The ring
$K_0^p(\Lambda)$ corresponds to this situation, since $K_0^p(\Lambda)=\zz\si
\oplus \zz\pr$. Moreover, $\zz\pr$ is a $\zz\si$-bimodule isomorphic to
$\zz\si$; we need to study the multiplicative bimodule structures of a ring
having underlying bimodule the ring itself.

\begin{pro}
\label{multiplic}
The multiplicative bimodule structures on a ring $A$ are in bijection with its
center $ZA$. The isomorphism classes of bimodule structures on $A$ are in
one-to-one correspondance with $ZA/U(ZA)$, where $U(ZA)$ is the group of units
of $A$.
\end{pro}
{\bf Proof:} Notice that an $A$-bimodule map $\varphi : A\otimes_A A\rightarrow
A$ is an  $A$-bimodule endomorphism of $A$  which  is completely determined by
the image of $1$, an element of the center of $A$. The associativity  always
holds since
$$\varphi\left(\varphi(a\oti b)\oti c\right)=abc\
\varphi(1)^2=
\varphi\left(a\oti\varphi(b\oti c)\right).$$
Finally an $A$-bimodule automorphism of  $A$ is given by a unit  of $ZA$, we
infer that two multiplicative bimodule structures $\varphi$ and $\varphi'$ are
isomorphic if and only if there exist a unit $u\in U(ZA)$ such that
$\varphi(1)=u\varphi'(1)$.

\vskip3mm
We define now a  canonical element $c$ of the integral group ring $\zz\si$ as
follows: let $P_1$ be the projective cover of the trivial module $k_\varepsilon$
(according to the preceding notation we have $P_1=P_{k_\varepsilon}$). The
element $c$ is provided by a composition series of $P_1$, more precisely
$c=\ov{P_1}$. Notice that for each $S\in\si$, we have $\ov{P_S}= Sc$ in the
group ring.

\begin{pro}
Let $\Lambda$ be a basic and split Hopf algebra, with structure group $\si$ (non
necessarily abelian) and canonical element $c$. Then $\kp$ is a ring isomorphic
to $\zz\si
\oplus
\zz\si$, where the first factor is the group ring and multiplication is given by :
$$(s_1,t_1)(s_2,t_2)=(s_1s_2, s_1t_2+t_1s_2+t_1ct_2).$$
\end{pro}

\begin{ex}
\rm
Let $q$ be a primitive $n$-th root of unity in $\cc$, and let ${{\sf
u}_q^+\!\left(sl_2\right)}$ be the corresponding finite dimensional half-quantum
group, which  structure group $\si$ is isomorphic to the cyclic group of order
$n$ (see \cite {hqg}). The element $c$ is the trace element in $\zz\si$, i.e.
the sum of all the group elements (the  structure of ${{\sf
u}_q^+\!\left(sl_2\right)}$-modules is given in \cite{aqqg}).
\end{ex}
\begin{cor} The canonical element
$c$ of a basic and split non semisimple Hopf algebra lies in the center of the integral group ring
of its structure group.
\end{cor}
In case the structure group is commutative there is a more explicit description
of the algebra ${\bf C}\otimes\kp$ . Recall that $\cc\si$ is a semisimple
commutative algebra isomorphic to $\cc^{\si}=\times_{S\in\si}\delta_S$, where
$\left\{\delta_S\right\}_{S\in\si}$ is the complete set of primitive orthogonal
idempotents provided by the Dirac masses on elements of $\si$. An explicit
isomorphism is provided by the Fourier transform once  generators of $\si$ and
corresponding roots of unity are given.
\vskip3mm
\noindent
{\bf Definition:} Let  $\Lambda$ be a basic and split Hopf algebra with abelian
structure group $\si$ and canonical element $c$.  The positive integer $r$ is
the number of non-zero coordinates of $c$ in the above decomposition of the
group algebra $\cc\si$. Equivalently, $r$ is the number of irreducible
characters of $\si$ which do not vanish at $c$.

\begin{th}
\label{struc}
Let $\Lambda$ be a basic and split Hopf algebra with abelian structure group of
order $s$. Then $\cc\ot\kp$ is isomorphic to the algebra $\cc^{2r}\times
\cc[\epsilon]^{s-r}$ where $\cc [\epsilon]$ is the algebra of dual numbers
 $\cc[X]/X^2$.
\end{th}
{\bf Proof:} We know that a multiplicative bimodule structure on a ring $A$
having underlying bimodule $A$ is determined by a central element up to
multiplication by a unit of the center. After normalization  we  have to
describe the algebra structure of $\cc^s\oplus\cc^s$ where the first factor
affords its natural algebra structure and the second is considered as a
$\cc^s$-bimodule with multiplicative structure provided by an element $c'$
having $r$ entries with value $1$, the others entries having zero value. For
simplicity we consider $\cc^s$ as the algebra of functions $\cc^E$ on the set
$E$, and $c'$ be given by a subset $F$ in $E$. A complete set of primitive
orthogonal idempotents is
$$\left\{\left(\delta_y,0\right)\right\}_{y\in E\setminus F}
\ \sqcup\
\left\{\left(\delta_x,-\delta_x\right),
\left(0,\delta_x\right)\right\}_{x\in F}.$$
Moreover, the elements $\left\{\left(0,\delta_y\right)\right\}_{y\in E\setminus
F}$ are orthogonal elements with zero square. Actually
$$\left(\delta_y,0\right)\left(0,\delta_y\right)=\left(0,\delta_y\right)\
\mbox{if}
\ y\in E\setminus F$$
$$\mbox{while}
\left(\delta_x,\delta_{-x}\right)\left(0,\delta_y\right)
=0=
\left(0,\delta_x\right)\left(0,\delta_y\right)
\ \mbox{for}\
y\in E\setminus F
\ \mbox{and}\
x\in F.$$ This formulas are equivalent to the structure result for $\kp$.

\begin{rem}
\rm
For ${{\sf u}_q^+\!\left(sl_2\right)}$, we notice that $c$  is the sum of all
the elements of the cyclic group of order $n$, hence  $r=1$. Then
$$\cc\ot K_0^p\left({\sf u}_q^+\!\left(sl_2\right)\right)
\ \cong\
\cc^{2} \times \cc[\epsilon]^{n-1}.$$
\end{rem}

Finally we  describe the virtual projective modules which have zero square.
Consider $G$ a finite abelian group with a given system of generators
$$G=<K_1,\dots,K_t\ \mid\ K_i^{n_i}=0\ \ K_iK_j=K_jK_i >,$$
and let $q_1,\dots,q_t$ be complex primitive roots of unity of order
$n_1,\dots,n_t$.  If $a=(a_1,\dots,a_t)$ is a $t$-uple of integers, we put
$q^a=q_1^{a_1}\cdots q_t^{a_t}$, and $K^a=K_1^{a_1}\cdots K_t^{a_t}$. Moreover
 $ab=\sum_{i=1}^{t}a_ib_i$ for  $b=(b_1\cdots b_t)$.

Consider the  non-degenerate bilinear form
$$\beta_G : {\bf C}G\times{\bf C}G\longrightarrow {\bf C}$$
given by the isomorphism $\alpha : {\bf C}G\rightarrow {\bf C}^G$ of Hopf
algebras, namely $$\beta_G(K^a,K^b)=q^{ab}.$$

Let  $\Lambda$ be a basic and split Hopf algebra with structure group $\si$ and
canonical element $c$. The ideal of projective modules in $\cc\ot\kp$ is
isomorphic to $(\cc\si)_c$ where the product is given by $t_1
.t_2
=t_1ct_2$. Let $B_c=\{K^x\in\si\ \mid\ \beta_{\si}(K^x,c)\neq 0\}$.

\begin{pro}
The nilradical of $\cc\ot\kp$ is the orthogonal vector space to $\cc B_c$ in
$\cc\si$.
\end{pro}

{\bf Proof :} Consider the corresponding product in $\cc^{\si}$,  given by
$f.g=f\alpha(c)g$. An element $f$ has zero square if $0=f.f=ff\alpha(c)$, which
only holds for $f\alpha(c)=0$. This means that for $K^x\in B_c$, we have
$f(K^x)=0$. The corresponding element in $\cc\si$ is orthogonal to each element
of $B_c$.

\begin{ex}
\rm
For $\uqpsl$ the element $c$ is the sum of all the elements of the cyclic group
of order $n$,  then $B_c=\{1\}$. The orthogonal space is the augmentation ideal
and any difference of two indecomposable projective modules has zero square.
\end{ex}

 \footnotesize

\vskip3mm
\noindent
D\'ep. de Math\'ematiques, Universit\'e de Montpellier 2, F-34095 Montpellier cedex 5.

\noindent
{\tt cibils@math.univ-montp2.fr}

\vskip2mm

\noindent
February 1998
\end{document}